\documentclass[11pt,bezier,amstex]{article}  % include bezier curves
\usepackage{color}              % Need the color package
\usepackage{epsfig}
\usepackage{graphicx}
\usepackage[]{amsmath}
\usepackage{amssymb}
\usepackage{epsfig}
\usepackage{hyperref}
\definecolor{MyDarkBlue}{rgb}{0,0.08,0.50}
\definecolor{BrickRed}{rgb}{0.65,0.08,0}

\hypersetup{
colorlinks=true,       % false: boxed links; true: colored links
    linkcolor=MyDarkBlue,          % color of internal links
    citecolor=BrickRed,        % color of links to bibliography
    filecolor=red,      % color of file links
    urlcolor=cyan           % color of external links
}

\newtheorem{Lemma}{Lemma}[section]
\newtheorem{Proposition}[Lemma]{Proposition}

\newtheorem{Theorem}[Lemma]{Theorem}

\newtheorem{Remark}[Lemma]{Remark}

\newtheorem{Construction}[Lemma]{Construction}
\newtheorem{Corollary}[Lemma]{Corollary}

\def\1{{\mathchoice {1\mskip-4mu\mathrm l}      % Blackboard bold 1
{1\mskip-4mu\mathrm l}
{1\mskip-4.5mu\mathrm l} {1\mskip-5mu\mathrm l}}}

\newcommand{\Qprob}{{\mathbb Q}}

\newcommand{\prob}{\mathbb{P}}

\newcommand{\FF}{\mathcal{F}}

\newcommand{\eps}{\varepsilon}

\newcommand{\qed}{\ \ \rule{1ex}{1ex}}

\newcommand{\Rbold}{{\mathbb{R}}}

\newcommand{\ind}[2]{1_{(e \in \pi(#1,#2))}}

\newcommand{\expec}{\mathbb{E}}

\def\ind{{\rm 1\hspace{-0.90ex}1}}

\newcommand{\eqn}[1]{\begin{equation} #1 \end{equation}}

\newcommand{\lbeq}[1]{\label{#1}}
\newcommand{\refeq}[1]{(\ref{#1})}
\newcommand{\sss}{\scriptscriptstyle}

\newcommand{\CE}{{\rm C}}

\newcommand{\op}{o_{\sss \prob}}
\newcommand{\Op}{O_{\sss \prob}}
\newcommand {\convd}{\stackrel{d}{\longrightarrow}}
\newcommand {\convp}{\stackrel{\sss {\mathbb P}}{\longrightarrow}}
\newcommand {\convas}{\stackrel{\sss a.s.}{\longrightarrow}}
\newcommand {\vep}{\varepsilon}

\numberwithin{equation}{section}
\newcommand{\BP}{{\rm SWT}}

\newcommand{\Bin}{{\rm Bin}}
\newcommand{\EXP}{{\rm Exp}}
\newcommand{\Poi}{{\rm Poi}}
\newcommand{\e}{{\rm e}}

\begin{document}
\author{Shankar Bhamidi
\thanks{ Department of Statistics and Operations Research,
304 Hanes Hall, University of North Carolina, Chapel Hill, NC. 27599, USA,
E-mail: {\tt bhamidi@email.unc.edu}} \and Remco van der Hofstad
\thanks{Department of Mathematics and
Computer Science, Eindhoven University of Technology, P.O.\ Box 513,
5600 MB Eindhoven, The Netherlands. E-mail: {\tt
rhofstad@win.tue.nl}}
\and
Gerard Hooghiemstra
\thanks{DIAM, Delft University of Technology, Mekelweg 4, 2628CD Delft, The
Netherlands. E-mail: {\tt g.hooghiemstra@tudelft.nl} } }

\title{First passage percolation on the Erd\H{o}s-R\'enyi random graph}

\maketitle

\begin{abstract}
In this paper we explore first passage percolation (FPP) on the
Erd\H{o}s-R\'enyi random graph $G_n(p_n)$, where each edge is
given an independent exponential edge weight with rate 1. In the
sparse regime, i.e., when $np_n\to \lambda>1,$ we find refined
asymptotics both for the minimal weight of the path between
uniformly chosen vertices in the giant component, as well as for the
hopcount (i.e., the number of edges) on this minimal weight path.
More precisely, we prove a central limit theorem for the
hopcount, with asymptotic mean and variance both equal to
$\lambda/(\lambda-1)\log{n}$. Furthermore, we prove that the
minimal weight centered by $\log{n}/(\lambda-1)$ converges in
distribution.

We also investigate the dense regime, where $np_n \to \infty$. We
find that although the base graph is a {\it ultra small}
(meaning that graph distances between uniformly chosen
vertices are $o(\log{n})$), attaching random edge weights changes
the geometry of the network completely. Indeed, the hopcount $H_n$
satisfies the universality property that whatever be the value of
$p_n$,  \ $H_n/\log{n}\to 1$ in probability and, more precisely,
$(H_n-\beta_n\log{n})/\sqrt{\log{n}}$, where
$\beta_n=\lambda_n/(\lambda_n-1)$, has a limiting standard normal
distribution. The constant $\beta_n$ can be replaced by 1
precisely when $\lambda_n\gg \sqrt{\log{n}}$, a case that has
appeared in the literature (under stronger conditions on
$\lambda_n$) in \cite{vcg-random-shanky, hofs-erdos-fpp}. We also
find bounds for the maximal weight and maximal hopcount between
vertices in the graph. This paper continues the investigation of
FPP initiated in \cite{vcg-random-shanky} and \cite{BHHS08}.
Compared to the setting on the configuration model studied in
\cite{BHHS08}, the proofs presented here are much simpler due to a
direct relation between FPP on the Erd\H{o}s-R\'enyi random graph
and thinned continuous-time branching processes.
\end{abstract}

\vspace{0.3in}


\noindent
{\bf MSC2000 subject classification.}
60C05, 05C80, 90B15.

\section{Introduction}
First passage percolation is one of the most fundamental problems
in probability theory. The basic motivation to study this
problem is the following. The goal is to model the flow of fluid
through some random medium. Suppose that we have a base graph on $n$
vertices which represents the available pathways for the fluid.
We attach to each edge in the graph some random edge weight, typically
assumed to be independent and identically distributed
(i.i.d.)~positive random variables with some probability density
function $f$. We then think of fluid percolating through the
network at rate $1$ from some source. Letting
$n\to\infty$, one is then interested in asymptotics of various
statistics of the flow through this medium. See e.g.,
\cite{hamm-welsh, howard} for results and a survey on the integer
lattice.

Our aim in this paper is to rigorously analyze first passage percolation (FPP) on the Erd\H{o}s-R\'enyi random graph
(ERRG) denoted by $G_n(p_n)$.
We shall see that for two randomly chosen vertices in the giant component
the hopcount, i.e., the number of edges
on the {\it shortest}-weight path between these vertices,
scales as $\log{n}$ and we shall find a central limit theorem (CLT) for this quantity.
We shall also find that the weight of the shortest-weight path re-centered by a constant
multiple of $\log{n}$ converges in distribution to some limiting random variable. We shall describe the
explicit distribution of the limit.
We shall also find lower bounds for the maximal optimal weight and hopcount
between vertices in the giant component.

In \cite{BHHS08}, we have investigated FPP on the configuration model (CM)
with degrees given by an i.i.d.~sequence with distribution function $F$, satisfying $F(x)=0,\, x<2$.
Consequently, all degrees are at least $2$ and the giant component contains $n-o(n)$ vertices, so that
with high probability ({\bf whp}), two uniformly chosen vertices are connected.
Furthermore, it was assumed that
for all $x\ge 0,$ there exist constants $c_1,c_2$, such that
    $$
    c_1 x^{1-\tau}\le 1-F(x)\le c_2 x^{1-\tau}, \quad \mbox{when} \quad \tau\in (2,3),
    $$
whereas for $\tau>3$, $F$ should satisfy, for all $x\ge 0$,
    $$
    1-F(x)\le cx^{1-\tau},
    $$
for some constant $c$.
Apart from self-loops and multiple edges, the CM, with a
binomial degree sequence, is not very different
from the ERRG. The main
challenge in studying the ERRG compared to the study of FPP on the CM in \cite{BHHS08} is threefold:
\begin{enumerate}
\item[(a)] The degree sequence in \cite{BHHS08} was assumed to
satisfy $F(x)=0,\, x<2$, so that all degrees are at least $2$,
with probability 1. As said this implies that uniformly chosen
vertices are {\bf whp} connected. In the present paper a uniformly
chosen pair of vertices does have a positive
probability of being connected, but this probability does not
equal $1$, so that we have to condition on the uniformly chosen
vertices to be in the giant component. Technically, this is a big
step forwards.

\item[(b)] We deal with the case where the average degree $\lambda_n=np_n\to\infty$, see
Corollary \ref{cor-lim} below, a scenario that is not contained in the sparse setting in \cite{BHHS08}.
This is the first time a result of this generality has been proved in the regime $np_n\to\infty$.

\item[(c)]
The  ERRG admits an elegant embedding in a
marked branching process in continuous time.
Consequently, the proofs are short and non-technical compared to those in \cite{BHHS08}. This
technique should prove to be useful in a number of other random graph models.
Moreover, we include a lower bound on the weight and length of the largest shortest-weight path for FPP on the ERRG.
\end{enumerate}

\section{Results}
\label{sec-res}
In this section we formulate our main results. Throughout this paper, we work on the
Erd\H{o}s-R\'enyi random graph (ERRG) $G_n(p_n)$, with vertex set $[n]=\{1, \ldots, n\}$
and edge set $E_n=\{(i,j):i,j\in [n], i\neq j\}$, and where every pair of vertices
$i\neq j$ is connected independently with probability $p_n$.
Furthermore, each edge $e\in E_n$ is equipped with an independent weight $E_e$, having an
exponential distribution with rate 1. We denote by $\convas$, $\convd$, and  $\convp$, convergence almost surely, in distribution, and
in probability, respectively. The symbols $o$, $O$ are the ordinary {\it Landau} symbols.
We say that a sequence of random variables $X_n$ satisfies $X_n=\op(b_n)$, $X_n=\Op(b_n)$, respectively,
if $X_n/b_n \convp 0$, $X_n/b_n$ is tight, respectively. We write that a sequence
of events $({\cal E}_n)_{n\geq 1}$ occurs \emph{with high probability} ({\bf whp})
when $\prob({\cal E}_n)=1-o(1).$ Binomial random variables are denoted by $\Bin(n,p)$,
where $n$ denotes the number of trials and $p$ the success probability.
We further notate by $\EXP(\mu)$, $\Poi(\lambda)$, respectively an exponentially distributed random variable with rate $\mu$,
and a Poisson random variable with mean $\lambda$.

In the Theorems \ref{theo:hopcount} and \ref{theo:weight} below, we shall investigate the
hopcount and weight of FPP on the ERRG.

\begin{Theorem}[CLT for hopcount]
\label{theo:hopcount}
Let $\lim_{n\to \infty } np_n= \lambda>1$ and define $\beta=\frac{\lambda}{\lambda-1}>1$.
Then, the hopcount $H_n$ between two uniformly chosen vertices,
conditioned on being connected, satisfies a central limit theorem
with the asymptotic mean and variance both equal to $\beta \log n$, i.e.,
\[
\frac{H_n -
\beta\log{n}}{\sqrt{\beta\log{n}}}
\convd Z,
\]
where $Z$ is a standard normal random variable.
\end{Theorem}

Now we consider the asymptotics of the minimal weight.

\begin{Theorem}[Limit distribution for minimal weight]
\label{theo:weight}
Let $\lim_{n\to \infty } np_n= \lambda>1$, and define $\gamma=\frac{1}{\lambda-1}$.
Then, there exists a non-degenerate real valued random variable
$X$ with distribution $\rho$, such that the minimal weight $W_n$ between two uniformly chosen vertices,
conditioned on being connected, satisfies
\[
W_n - \gamma \log{n}\stackrel{d}{\longrightarrow} X.
\]
\end{Theorem}

\begin{Remark}[Joint convergence]
The proof shall reveal that the convergence in Theorems \ref{theo:hopcount} and \ref{theo:weight}
holds \emph{jointly}, with the limits being independent.
\end{Remark}

Let us identify the limiting distribution $\rho$.
Consider a continuous-time Galton-Watson branching-process
with $\Poi(\lambda)$ offspring, where the individuals have an exponential life time
with rate 1. Denote the number of alive individuals at time $t$ by
$N(t)$. It is well known \cite{athreya} that there exists $\alpha>0$
such that $\e^{-\alpha t} N(t)$ has an almost sure limit $W$. For the case under
consideration it is readily verified that the Malthusian parameter $\alpha=\lambda-1>0$.
Let $D$ denote a $\Poi(\lambda)$ random variable representing the number of offspring
from one individual, and, conditioned on $D$, denote by $E_1,E_2,\ldots,E_D$, the exponential lifetimes
of the offspring. Then the limit $W$ satisfies the stochastic equation
\eqn{
\label{stoch-eq}
W\stackrel{d}{=}\sum_{i=1}^D \e^{-\alpha E_i} W_i,
}
where on the right-hand side all the involved random variables are independent
with $W_i$ being copies of $W$. Hence taking conditional expectations w.r.t. $D$, it is
seen that
$\phi(t) = \expec[\e^{-tW}]$ satisfies the functional relation
\eqn{
\label{func-eq}
\phi(t) = \expec
\Big[
\Big(
\expec[\e^{-t \e^{-\alpha E_1} W}]
\Big)^D
\Big]=
\exp\left(-\lambda\int_0^\infty \left[1-\phi(t\e^{-(\lambda-1)x})\right]\e^{-x} dx\right).
}
From this, it is quite easy to prove that the random variable $W$ has an atom at zero
of size $p_\lambda$, where $p_{\lambda}$ is the smallest non-negative solution of the equation
\[
p_\lambda = \exp(-\lambda(1-p_\lambda)),
\]
i.e., $p_{\lambda}$ is the extinction probability of a branching process with
$\Poi(\lambda)$ offspring.
Furthermore, the random variable $W$, conditioned to be positive, admits a continuous density on $\Rbold^+$, and we denote by
\eqn{
W_\lambda \stackrel{d}{=} (W|W>0).
\label{eqn:wl}
}

To construct $\rho$ we shall need the following random variables:

(a) Let $W_\lambda^{\sss (1)}, W_\lambda^{\sss (2)}$ be independent and identically distributed as
${W_\lambda}$,
\\
(b) Let $E \stackrel{d}{=}\EXP(1)$ be independent of $W_\lambda^{\sss (i)}$.

In terms of these random variables, the random variable $X$ with distribution $\rho$ satisfies
\[
X \stackrel{d}{=} - \gamma\log{(\gamma W_\lambda^{\sss (1)})}-\gamma\log{(\gamma W_\lambda^{\sss (2)})}
+\gamma \log(E).
 \]

We next study the \emph{dense graph setting,} where $np_n=\lambda_n\rightarrow \infty$.
The proof in this setting follows the same lines as that of Theorems
\ref{theo:hopcount}-\ref{theo:weight}, and therefore we will only give a sketch of proof.
Observe that for $np_n\to \infty$, any pair of vertices is {\bf whp} connected.
It is not hard to see that, in this case, the giant component consists of $n(1-o(1))$
vertices, and that the graph distance between two uniformly chosen vertices
is $\log{n}/\log{\lambda_n}=o(\log{n})$, so that the random graph is \emph{ultra small}.
In the statement of the result, we denote by $H_n$ and $W_n$,
respectively, the hopcount and minimal weight of the shortest-weight path between two uniformly chosen vertices.

\begin{Corollary}[Limit of hopcount and weight for $np_n\to\infty$]
\label{cor-lim}
Set $\lambda_n=np_n$. For $\lambda_n\to\infty$, as $n\to \infty$,
\\(a)
\eqn{
\lbeq{hopconlamn}
\frac{H_n -\beta_n \log{n}}{\sqrt{\log{n}}} \convd Z,
}
 where $\beta_n=\lambda_n/(\lambda_n-1)$ and $Z$ is a standard normal random variable;
\\(b)
\eqn{
\lbeq{weightconlamn}
(\lambda_n-1) W_n-\log n\convd {\tilde X},
}
where the random variable ${\tilde X}\stackrel{d}{=}M_1+M_2-M_3$, with $M_1$, $M_2,M_3$  independent
Gumbel random variables, i.e.,
$\prob(M_i\le x)=\Lambda(x)=\exp(-e^{-x}),\, 1\le i \le 3,$ for all $x\in \Rbold$.
\end{Corollary}
\begin{Remark}[Dense setting]
(a) In Part (a) of Corollary \ref {cor-lim} the centering $\beta_n \log{n}$ can be replaced by $\log{n}$ if and only if
$\lambda_n/\sqrt{\log{n}}\to \infty$. In this case, the hopcount has the same
limiting distribution as on the complete graph (see e.g., \cite{hofs-erdos-fpp}).\\
(b) Note that the distribution of $H_n$ for the case where $np_n/(\log{n})^3 \to \infty$ was obtained in
\cite{hofs-erdos-fpp}. In \cite[Theorem 5]{vcg-random-shanky} limit results for $H_n/\log{n}$ and
$np_n W_n-\log{n}$ were obtained in the case where $\lim \inf_{n\to \infty }(np_n)/(\log{n})=a$, with
$1<a\le \infty$. Thus, the present paper essentially completes the study of FPP on ERRGs
in all regimes.
\end{Remark}

Theorems \ref{theo:hopcount}-\ref{theo:weight} state that for
uniformly chosen connected vertices, the weight and length of the
optimal path between them scales as $\gamma \log{n}$ and
$\beta\log{n}$, respectively. The following shows the existence of
pairs of vertices with a much larger optimal path
weight and hopcount:

%%%%%%%%%%%%%%%%%%%%%%%%%%%%%%%%%%%%%%%%%%%%%%%%%%%%%%%%%%%%%%%%%%%%%%%%%%%%%%%%%%%
\begin{Theorem}[Existence of long paths]
\label{theo:extrema}
Let $np_n\to \lambda>1$ and define by  $c(\lambda), d(\lambda)$ the constants
    \eqn{
    c(\lambda) = \frac{1}{\lambda-1} + \frac{2}{\log{|\mu_\lambda|}},
    \qquad
    d(\lambda)=\frac{\lambda}{\lambda-1} + \frac{2}{\log{|\mu_\lambda|}},
    }
where $\mu_\lambda$ is the dual of $\lambda$, i.e., the unique $\mu\in (0,1)$ such that
    \eqn{
    \lbeq{dual-def}
    \mu\e^{-\mu}=\lambda\e^{-\lambda}.
    }
Then for any given $\eps$, with high probability, there exists a pair of vertices say $i^*, j^*$ in the giant component
such that the weight of the optimal path $W_n(i^*, j^*)$ and its number of edges $H_n(i^*,j^*)$ satisfy the inequalities
    \eqn{
    W_n(i^*,j^*) \geq (1-\eps)c(\lambda) \log{n},
    \qquad
    H_n(i^*,j^*) \geq (1-\eps)d(\lambda) \log{n}.
    }
\end{Theorem}
We conjecture that  the above result is optimal, i.e., you can get between any pair of
vertices within weight $c(\lambda)\log{n}$ and with at most $d(\lambda)\log{n}$ hops.

In \cite{BHHS08}, we have proven results parallel to those in
Theorem \ref{theo:hopcount} and \ref{theo:weight} for the
configuration model (CM) with degrees given by an i.i.d.~sequence
with distribution function $F$, where $F(x)\le c x^{1-\tau}, \,
\tau>3$, for some constant $c$ and all $x\ge 0$. We found that the asymptotic hopcount between two
uniformly chosen vertices converges to a normal distribution with
mean and variance equal to $\frac{\nu}{\nu-1} \log n$, where $\nu
=\expec[D(D-1)]/\expec[D]$ and where the random variable $D$ is
distributed as the degree distribution $F$. Note that in Theorem
\ref{theo:hopcount} the role of $\nu$ is taken over by the
parameter $\lambda$. This is not surprising, since for a
$\Poi(\lambda)$ variable $D$, we have $\nu
=\expec[D(D-1)]/\expec[D]=\lambda$.

\section{Setting the stage for the proofs}
\label{sec-proof}
A rough idea of the proof is as follows. Fix two vertices, say $1$ and $2$, in the giant component.
Think of fluid emanating from these two sources {\it simultaneously} at rate $1$,
so that at time $t$, $\FF^{\sss (i)}(t)$ is the flow cluster from vertex $i$ and includes the
minimal weight paths to all vertices wetted at or before $t$ from vertex $i,\, i=1,2$.
When these two flows collide (namely, when the flow from one of the sources reaches the other flow)
via the formation of an edge $(v_1, v_2)$ between two vertices $v_1\in \FF^{\sss (1)}(\cdot)$
and $v_2\in \FF^{\sss (2)}(\cdot)$, then the shortest-weight path between the two vertices has been
found. This collision time, which we denote by $S_{12}$,  tells us that the weight between
the two vertices, $W_n(1,2)$, equals
    \eqn{
    W_n(1,2) = 2S_{12}.
    \label{eqn:wn}
    }
Furthermore, if $G_n(v_i), i=1,2 $, denotes the number of edges between the source $i$
and the vertex $v_i$ along the tree $\FF^{\sss(i)}(S_{12})$, then the hopcount $H_n(1,2)$ is given by
    \eqn{H_n(1,2) = G_n(v_1)+ G_n(v_2) +1.
    \label{eqn:hn}
    }

The above idea is indeed a very rough sketch of our proof. In the paper we embed the flow on the ERRG in a
continuous-time marked branching process (CTMBP), where the offspring distribution is
binomial with parameters $n-1$ and $p$ (see Section \ref{CTBP}).
With high probability, the marks in the CTMBP correspond to the vertices in the ERRG.
We denote by $\{\BP^{\sss (i)}_m\}_{m\ge 0}$ the marks of the individuals wetted by the flow after $m$ splits
and the arrival times of these splits, where the superscript $i,\, i=1,2,$ denotes the root $i$ of the flow.
It is not mandatory to let the flows grow simultaneously, and for technical reasons, the proof is simpler if we first grow
$\BP_{m}^{\sss(1)}$ to a size $a_n=\lceil \sqrt{n}\rceil$. After this, we grow $\BP_{m}^{\sss(2)}$, and we stop as soon as a mark of
$\BP_{a_n}^{\sss(1)}$ appears in $\{\BP_m^{\sss(2)}\}_{m=0}^{\infty}$. The size $a_n=\lceil \sqrt{n}\rceil$ is the correct one, since
if both flows are of size approximately $\sqrt{n}$, the probability that the second flow finds a mark of the first flow
in some time interval of positive length is of order $1$.

The CLT for the distance $G_m$ between the $m^{\rm th}$-wetted mark and the root as well as the limit
distribution for the weight of the path between this mark and the root are given in Section \ref{HWRNV}.
This theory is based on the asymptotics for Bellman-Harris processes which will be developed in
Section  \ref{sec:bh}.

In Section \ref{TCT} we investigate the connection time $\CE_n$, i.e., the random time until the second flow starting from
mark 2 hits the flow of size $a_n$, which started from mark 1. We prove that $\CE_n/a_n$ converges in distribution to an 
$\EXP(1)$ random variable. This can informally be understood as follows.
The number of distinct marks in $\BP_{a_n}^{\sss(1)}$ is $a_n(1+o(1))$. Each of these marks
is chosen with probability $1/n$ in $\BP_{m}^{\sss(2)}$, so that the first time that
any of these marks is chosen is close to a geometric random variable with
success parameter $a_n/n$. This random variable is close to $n/a_n \EXP(1)$,
so that indeed $\CE_n/a_n\convd \EXP(1)$.

We furthermore show that, conditioned on $\CE_n$, the hopcount $H_n$ is, {\bf whp},
the independent sum of $G^{\sss (1)}_{U_n,a_n}$ and $G^{\sss (2)}_{\CE_n}$. Here
$G^{\sss (1)}_{U_n,a_n}$ denotes the distance between root 1 and a mark in $\BP_{a_n}^{\sss(1)}$
that is chosen uniformly at random.

With the above sketch of proof in mind, the remainder of the paper is organized as follows:
\begin{itemize}
\item[(1)] In Section \ref{sec:bh} we will analyze various properties of a Bellman-Harris process
conditioned on non-extinction, including times to grow to a particular size and the
generation of individuals at this time. In this continuous-time branching process the offspring will have a Poisson distribution.

\item[(2)] In Section \ref{BPTIE}, %carefully
we introduce marked branching process trees with binomially distributed offspring and make the connection between these trees
and the ERRG, by thinning the marked branching process tree.

\item[(3)] In Section \ref{CTBP},
we replace the general weights on the edges by exponential ones so that we
end up with a {\it continuous-time} marked branching processes (CTMBP) with binomial offspring.
We further focus on some of the nice properties of the involved random variables, that are a consequence of the
memoryless property of the exponential distribution.

\item[(4)] By coupling the CTMBP with binomial offspring to a  Bellman-Harris process with Poisson offspring, we
deduce in Section \ref{HWRNV} the limit theorem for the generation $G_m$ and
the weight $A_m$ from the results in Section \ref{sec:bh}.

\item[(5)] In Section \ref{TCT}, we present a refined analysis of the connection time
of the two flow clusters.

\item[(6)] Finally, in Section \ref{sec-proof-main-res}, we complete the proofs of our main results.

\end{itemize}
The idea of the argument is quite simple but making these ideas rigorous
takes some technical work because of the issue of conditioning on being in
the giant component.

\subsection{Asymptotics for Bellman-Harris processes}
\label{sec:bh}
Here we shall construct random trees where each vertex lives for an exponential amount of time, dies
and gives birth to some number of children.
Fix a sequence of non-negative integers $d_1, d_2, \ldots$. For future reference, let
$s_i = \sum_{j=1}^i d_j - (i-1)$. We shall make the following blanket assumption
\eqn{
s_i > 0,\quad \mbox{ for all } i\geq 1.
 \label{eqn:posts}
 }
The condition in \refeq{eqn:posts} is equivalent to the fact
that the tree where the $i^{\rm th}$ vertex has degree
$d_i$ is infinitely large. Consider the following
{\it continuous-time} construction of a random tree:
\begin{Construction}[FPP on a tree]
\label{constr:cont-time}
{\it
\begin{enumerate}
\item[(1)] Start with the root which dies immediately giving rise to
$d_1$ alive offspring;
\item[(2)] each alive offspring lives for an $\EXP(1)$ amount of time, independent of all other
randomness involved;
\item[(3)] when the $i^{\rm th}$ vertex dies it leaves behind $d_i$ alive offspring.
\end{enumerate}
}
\end{Construction}

In terms of the above construction, we can identify $s_i$ as the number of alive
vertices after the $i^{\rm th}$ death. Let $T_1,T_2,\ldots,T_m$ be the time
spacings between the $(i-1)^{\rm st}$ and $i^{\rm th}$ death $1\le i \le m$. We shall
often refer to $A_i=T_1+T_2+\ldots+T_i,\,1\le i \le m,$ as the time of the $i^{\rm th}$ split.
Note that at the graph topology level, the above procedure is the same as the
following construction, where exponential life times do not appear:

\begin{Construction}[Discrete-time reformulation of FPP on a tree]
\label{constr:det}
The shortest-weight \- graph on a tree with degrees $\{d_i\}_{i=1}^{\infty}$
is obtained as follows:
\begin{enumerate}
\item[(1)] At time $0$, start with one alive vertex (the initial ancestor);
\item[(2)] at each time step $i$, pick one of the alive vertices at random,
this vertex dies giving birth to $d_i$ children.
\end{enumerate}
\end{Construction}
Let $G_m$ denote the graph distance between the root and a uniformly chosen vertex among all alive vertices at step $m$.
We quote the following fundamental result from
\cite{buhler}.
\footnote{A new {\it probabilistic} proof is given in \cite{BHHS08} ,
since there is some confusion in  comparing the
definition $s_i$ in this paper and the definition of $s_i$
given in \cite[below (3.1)]{buhler}.}

\begin{Proposition}[Shortest-weight paths on a tree]
\label{prop:gen}
Pick an alive vertex at time $m\geq 1$ uniformly at random from all vertices alive at this
time. Then\\
(a) the generation of the $m^{\rm th}$ chosen vertex is equal in distribution to
    \eqn{
    \lbeq{Gm-def}
    G_m = \sum_{i=1}^m I_i,
    }
where $\{I_i\}_{i=1}^{\infty}$ are independent Bernoulli random variables with
    \eqn{
    \lbeq{bi-def}
    \prob(I_i = 1) = d_i/s_i,\qquad 1\le i\le m.
    }
(b) the weight of the shortest-weight path between the root
of the tree and the $m^{\rm th}$ chosen  vertex is  equal
in distribution to
    \eqn{
    \lbeq{Tm-def}
    A_m=T_1+T_2+\ldots+T_m \stackrel{d}{=} \sum_{i=1}^m E_i/s_i,
    }
where $\{T_i\}_{i=1}^{\infty}$ are i.i.d.\ exponential random
variables with rate $s_i$ and hence $T_i$ is equal in distribution to $E_i/s_i$.
\end{Proposition}

\medskip
The above proposition states that the random variable $G_m$ is the sum of $m$ independent Bernoulli
random variables. It is known that in a wide variety of settings, such quantities
essentially follow a CLT.
We shall be interested in showing that the standardization $\tilde{G}_m$ of $G_m$
converges in distribution to a normal random variable $Z$, and the way we shall prove this is
by showing that the Wasserstein distance between the distribution $\tilde{G}_m$
and that of $Z$ goes to $0$.
Denote by  $\mbox{Wass}$ the metric on the space of probability measures
defined as
    \[
    \mbox{Wass}(\mu_1,\mu_2) = \sup\left\{\left|\int_\Rbold g d\mu_1 - \int_{\Rbold} g d\mu_2\right|:
    g \mbox{ is $1-$Lipschitz and bounded}  \right\}.
    \]
It is well known \cite[Theorem 11.3.3]{dudley} that this metric on the space of probability measures
induces the same topology as the topology of weak convergence.
Thus showing that $\tilde{G}_m\stackrel{d}{\longrightarrow} Z$
is equivalent to showing that $\mbox{Wass}(\mu_m, \mu) \to 0$ where $\tilde{G}_m \sim \mu_m$
and $Z\sim \mu$; here $X\sim \nu$ means that the law of the random variable $X$ is equal to $\nu$.

The following result is simple to prove using any of the methods
for showing CLTs for sums of independent random variables,
e.g. Stein's method \cite{stein-chen}.
\begin{Lemma}[CLT of hopcount in tree]
\label{lemma:normal}
As before let $ G_m$ denote the generation to a uniformly chosen
vertex alive at time $m$.  Define the sequence of random variables
    \eqn{
    \tilde{G}_m = \frac{G_m - \sum_{i=1}^m \rho_i}
    {\sqrt{\sum_{i=1}^m \rho_i(1-\rho_i)}},
    }
where $\rho_i = d_i/s_i$.
Let $\mu_m$ denote the distribution of ${\tilde G}_m$ and $\mu$
the distribution of a standard normal random variable. Then
    \eqn{
    \mbox{Wass}(\mu_m,\mu)
    \leq \frac{3}{\sqrt{\sum_{i=1}^m \rho_i(1-\rho_i)}}.
    }
\end{Lemma}

{\bf Proof.}
Combining \cite[Theorem 3.1 and 3.2]{stein-chen} yields
    \eqn{
    \label{wass-bnd}
    \mbox{Wass}(\mu_m,\mu)\leq 3\sum_{i=1}^m \expec[|\eta_i|^3],
    }
where
    \[
    \eta_i=\frac{I_i-\rho_i}{\sqrt{\sum_{i=1}^m \rho_i(1-\rho_i)}},
    \]
and where  $I_i$ is a $\Bin(1,\rho_i)$ random variable. From this we obtain
    \[
    \sum_{i=1}^m \expec[|\eta_i|^3]=\frac1{[\sum_{i=1}^m \rho_i(1-\rho_i)]^{3/2}} \sum_{i=1}^m \rho_i(1-\rho_i)\{\rho_i^2+(1-\rho_i)^2\}
    \le
    \frac{1}{\sqrt{\sum_{i=1}^m \rho_i(1-\rho_i)}},
    \]
using $\rho_i^2+(1-\rho_i)^2\le 1$. \hfill\qed

\smallskip

From now on we take $D_1,D_2,\ldots$ i.i.d.~where $D_j$ has a $\Poi(\lambda)$ distribution with $\lambda>1$.
Furthermore, we consider the continuous time construction of a random tree as in
Construction \ref{constr:cont-time},
with $d_i=D_i$ and $S_i = \sum_{j=1}^i D_j - (i-1)$.
Since the $D_j $ are now \emph{random}, for any fixed time $m\geq 1$, there is now some
non-zero probability that $S_j = 0$ for some $1\leq  j\leq m$. However, note that for any $m$,
    \[
    \prob\left(S_i > 0 \hspace{4mm} \forall i = 1,2,\ldots, m\right) \geq 1-p_\lambda,
    \]
where $p_\lambda$ is the extinction probability of a Galton-Watson branching process
with $\Poi(\lambda)$ offspring.
Let $G_m$ be the generation of a randomly chosen alive individual among all alive
vertices at time $m$, conditional on $S_i > 0$ for all $1\leq i\leq m$.
Let $A_m$ be the time for the $m^{\rm th} $ split to happen in the
continuous-time construction (see Construction \ref{constr:cont-time}).
Then we have the following asymptotics:

\begin{Theorem}[CLT for hopcount on trees conditioned to survive]
\label{theo:poisson}
Conditioned on $S_i > 0$ for all $1\leq i\leq m$, the following asymptotics for $m\to \infty$ hold:\\
(a) the generation $G_m$ satisfies a CLT, i.e.,
    \begin{equation}
    \label{weakconvG}
    \frac{G_m -  \beta\log{m}}
    { \sqrt{\beta\log{m}}} \convd Z,
    \end{equation}
where $Z$ is standard normal and $\beta=\frac{\lambda}{\lambda -1}$;
\\(b) the random variable $A_m$
satisfies the asymptotics:
\begin{equation}
\label{weakconvT}
A_m - \gamma\log{m}\convas
-\gamma\log{\gamma W_\lambda},
\end{equation}
where $W_\lambda$ has a distribution given by \eqref{eqn:wl}.
The limits in (a) and (b) also hold jointly, where the limits are independent.
\end{Theorem}

The issue of conditioning is a technical annoyance as it removes the
independence of increments of the random walk $\{S_i\}_{i=1}^m$. However what should be intuitively clear is that, for
a random walk with positive drift,
conditioned to be positive is, {\bf whp}, the same as conditioning the path to be positive in the first $w_m$ steps,
where we assume $w_m = o(\log\log{m})$. Indeed, by the time $w_m$,
the walk has reached height approximately $ (\lambda-1) w_m$ and if we just add independent
Poisson random variables from this stage onwards, then the random walk will, {\bf whp}, remain
positive, since the walk has reached a high level by this time, and the
probability that it would reach $0$ by time $m$ is exponentially small in $w_m$.
This idea is made precise in the following construction.
First we shall need some notation. Fix $w_m =o( \log\log{m}) \to \infty$ and let $Y_i = D_i-1$ where the $D_i$ are i.i.d.~$\Poi(\lambda)$ random variables.
Let $S_j =1+\sum_{i=1}^j Y_i$.

Now consider the following construction
\begin{Construction}
\label{construction-delay}
Take $S_{1}^{*}, S_{1}^{*}, \ldots, S_{w_m}^{*}$
equal in distribution to $S_1,S_2,\ldots,S_{w_m}$ conditioned on being positive.
Furthermore, let $Y_{w_m+1}, Y_{w_m+2}, \ldots Y_{m}$ be i.i.d.~distributed as
$\Poi(\lambda)-1,$ independent of $S_{1}^{*}, S_{1}^{*}, \ldots, S_{w_m}^{*}$ .
Define a new sequence $S_1, \ldots, S_m$ by
    \eqn{
     S_i
     =\left\{
     \begin{array}{ll}
     S_{i}^{*}, & i\leq w_m\\
     S_{w_m}^{*}+\sum_{j=w_m+1}^i Y_j,& w_m+1\le i\le m.
     \end{array}
     \right.
    }
\end{Construction}
Then the following proposition yields good bounds
on the behavior of this random walk:

\begin{Proposition}[Good bounds for the conditioned random walk]
\label{prop:error-bds}
The sequence in Construction \ref {construction-delay} satisfies the following regularity properties:
\\(a) $\prob(S_1>0,\ldots,S_m>0)=1-o(1)$, as $m\to \infty$, and there exists a coupling such that
$\prob(S_j={\tilde S_j},\, 1\le j \le m)=1-o(1)$, as $m\to \infty$, where ${\tilde S}_j, 1\le j \le m,$
is equal in distribution to the random walk $R_j=1+\sum_{i=1}^j (D_i-1), 1\le j \le m$,
conditioned on $R_1>0,\ldots,R_m>0$.
\\(b) There exists a constant $C>0$ such that, {\bf whp}, for all $j> \log\log{m}$,
    \eqn{
    (\lambda-1)j - C\sqrt{j\log{j}}
    \leq S_j \leq (\lambda-1)j +  C\sqrt{j\log{j}}.
    }
\\(c) Let $\rho_i = D_i/S_i$. Then
    \eqn{
    \sum_{i=1}^m \rho_i(1-\rho_i)= \beta\log{m} + \Op(\log{\log{m}}).%, \quad \sum_{i=1}^m \rho_i^2  = \Op(\log{\log{m}}).\\
    }
\end{Proposition}
Assuming Proposition \ref{prop:error-bds}, let us show how to prove Theorem \ref{theo:poisson}:

{\bf Proof of Theorem \ref{theo:poisson}.}
In order to prove Part (a), we note that according to Lemma \ref{lemma:normal} it  suffices to show
that $\sum_{i=1}^m \rho_i(1-\rho_i) \to \infty$, as $m\to \infty$. This is immediate from Part (a) and (c) of Proposition \ref{prop:error-bds}.

The almost sure convergence of $A_m=T_1+\ldots+T_m$ in \eqref{weakconvT} is a little more difficult.
We refer to \cite[Theorem 2, p.~120]{athreya} where it has been
proved that the time of the $m^{\rm th}$ split,  $A_m$, when the branching process has $\Poi(\lambda)$
offspring satisfies the asymptotics
    \eqn{
    \label{asymp-split}
    N(A_m) \e^{-(\lambda-1)A_m} \convas W,
    }
where $N(t)$ denotes the number of alive individuals at time $t$.
Observe that \refeq{asymp-split} can easily be deduced from the convergence
$\e^{-(\lambda-1)t}N(t)\convas W$, by substituting $A_m$ for $t$, and using that
$A_m\convas \infty$. Conditioning on the random walk to stay positive yields:
\eqn{
%\label{asymp-split}
\left(N(A_m) \e^{-(\lambda-1)A_m}|S_1>0,\ldots,S_m>0\right) \convas W_{\lambda}.
}
Now $N(t)$ is the number of {\it alive} individuals at time $t$, so after $m$ splits
$N(A_m)/m=S_m/m \convas \lambda-1$, by the strong law.
Hence, with $\gamma=1/(\lambda-1),$
    \eqn{
    \big(A_m| S_1>0,\ldots,S_m>0\big) -\gamma \log m \convas -\gamma\log{\gamma W_\lambda}.
    \label{eqn:cond-lt}
    }
\vskip-0.5cm

\hfill\qed

{\bf Proof of Proposition \ref{prop:error-bds}.} We start with the proofs of Part (b) and Part (c).
Note that $w_m=o(\log \log m)$, and that the random variables $D_j$ with $j>w_m$ are independent.
The proof of (b) is then straightforward from the large deviation properties of the Poisson
distribution and a simple union bound.
For (c) we write
    \[
   \sum_{i=1}^m \rho_i= \sum_{i=1}^m \frac{D_i}{S_i}=
    \sum_{i\le \log\log m} \frac{D_i}{S_i}+\sum_{i=\lceil \log\log m \rceil}^m \frac{D_i}{S_i}.
    \]
For the second sum on the right-side we combine the bounds in (b) and use Chebyshev's
inequality on the i.i.d.~Poisson random variables $D_i,\, i\ge \log\log m$, to obtain
    \[
    \sum_{i=\lceil \log\log m \rceil}^m \frac{D_i}{(\lambda-1)i} =\beta \log m+\Op(1).
    \]
The first sum on the right-side is obviously of order $\Op(\log\log m)$. This proves (c).

For the first statement in (a), we observe that $\{S_i>0~\forall i=1, \ldots, w_m\}$
and $\{S_m\geq x\}$ are both increasing in the i.i.d.~random variables $\{D_i\}_{i=1}^m$.
Therefore, by the FKG-inequality \cite{ForKasGin71} and for any $x\geq 0$,
we obtain that $\prob(S_m\geq x\mid S_i>0 \, \forall i=1, \ldots, w_m)\geq \prob(S_m\geq x).$
Therefore,
    \eqn{
    \prob(S_m\geq l_m\mid S_i>0 \, \forall i=1, \ldots, w_m)\geq \prob(S_m\geq l_m)=1-o(1),
    }
when we take $l_m =(\lambda-1) w_m -  C\sqrt{w_m\log{m}}$, and use Part (b).

We now turn to the coupling statement in Part (a). Denote by $\prob_l$ the probability distribution of the random walk starting at $l$. Then, using again
the notation $\{R_j\}_{j\ge 0}$ for an unconstrained random walk with step size $Y_k$,
\begin{eqnarray}
\label{martingalebnd}
\prob_{l_m}(\min_{0\le j \le m} S_j \leq 0)
\le \prob_{l_m}(\min_{ j \ge 0} R_j \leq 0)=\prob(\max_{j\ge 0} {\tilde R}_j\le l_m)\le \e^{-\theta^* l_m}=o(1),
\end{eqnarray}
where ${\tilde R}_j=-R_j, \, j\ge 0,$ and where $\theta^*>0$ is the unique {\it positive} solution
of the equation
\begin{equation}
\label{ross}
\theta-\lambda+\lambda \e^{-\theta}=0.
\end{equation}
Indeed, the final inequality in \refeq{martingalebnd} is obtained using the martingale 
$Z_n=\prod_{i\le n} \e^{-\theta^*Y_i} $, which
is by \refeq{ross} the product of independent unit mean random variables \cite[p.\ 342]{rossboek}.

Finally, if we take ${\tilde S}_j(\omega)=S^*_j(\omega)$ for $j<w_m$ and all $\omega$
in the probability space, whereas we take for ${\tilde S}_j(\omega)=S^*_{w_m}(\omega)+\sum_{i=w_m+1}^j Y_i(\omega)$, for
$j>w_m,$ and those $\omega$ with $l_m+\sum_{i=w_m+1}^j Y_i(\omega)>0, \, w_m<j\leq m$ and ${\tilde S}_j(\omega)=0$ otherwise,
we obtain our desired coupling. This completes the proof of Part (a) and, since the proof of Part (b) and Part (c)
has been given above, we are done. \hfill \qed

\section{Branching process trees and imbedding of ERRG}
\label{sec-BPembed}

We introduce marked branching process trees in Section \ref{BPTIE}. In Section \ref{CTBP}, we
attach the exponentially distributed weights to the edges and obtain in this way a continuous-time marked branching
process (CTMBP) with binomially distributed offspring.
Section \ref{HWRNV} describes the distribution of the height and weight of
minimal weight paths in these CTMBP. In Section \ref{TCT}, we treat the
connection between the CTMBP and the ERRG.

\subsection{Marked branching process trees}
\label{BPTIE}
In this section we describe a  marked branching process tree (MBPT) with a binomially distributed offspring, and we show
how this MBPT should be thinned to obtain the connection with the ERRG. We then attach i.i.d.~weights having probability density $f$ on $(0,\infty)$
to both the branching process tree and the
ERRG. We interpret these weights as distances between marks, vertices, respectively.
We close the section with a proof that for each time $t\ge 0,$ the set of marks that can be reached
in the thinned process within time $t$ from the root is identical in distribution to
the set of vertices that can be reached within time $t$ from the initial vertex $i_0$ in the ERRG.

The set of  marks is denoted by $[n]=\{1,2,\ldots,n\}$. In the branching process tree, we start with a single individual with mark
$i_0\in [n]$, and we put $I_0=\{i_0\}$. This individual reproduces a binomially distributed number of offspring $X_1$, with parameters $n-1$ and $p$.
We attach each of these $X_1$ children to
their father by a single edge, so that we obtain a tree structure. Each edge is assigned a weight taken
from an i.i.d.~sample having probability density $f$, with support $(0,\infty)$. To obtain the marks
of the offspring we take at random a sample of size $X_1$ from the set $[n]\setminus \{i_0\}$. The corresponding
set of marks is denoted by $I_1$.
Hence the marks are chosen uniformly at random from a set of size  $n-1$ and are all different.
Moreover the mark of the father ($i_0$) is not present among the set of marks $I_1$ of the offspring.
Observe that the marks can be seen as vertex numbers of the ERRG; from this viewpoint, the children of
individual $i_0$, given by the vertex-set $I_1$,  are the direct neighbors of the vertex $i_0$ in the ERRG.

We proceed by taking the individual for which the corresponding edge-weight is minimal. In our language this individual is
reached after $T_1$ time units, where $T_1$ is equal to the length of the minimal edge weight, since the  fluid
percolates at rate $1$. The mark reached by time $T_1$ is denoted by $i_1$. The reproduction process for $i_1$ is identical to that of $i_0$: we take
a $\Bin(n-1,p)$ number of children (denoted by $X_2$) and connect each of them by a single edge to $i_1$. The $X_2$ `new' edges are supplied by weights taken i.i.d.~with
probability density $f$. We update the `old' edges by subtracting $T_1$ from their respective weights, because $T_1$
time units have been used. We then complete this step by assigning marks to the new individuals,
taken uniformly at random from the set $[n]\setminus \{i_1\}$, this mark set
is denoted by $I_2$.

We proceed by induction: suppose that we have reached the individuals
$i_0,i_1,\ldots,i_k$, in this order, so that the (updated) weight attached to the edge with endpoint $i_k$
was the minimal edge weight at time  $A_k=T_1+T_2+\ldots+T_k$. We now form
$X_{k+1}\stackrel{d}{=}\Bin(n-1,p)$ new edges and attach them to $i_k$. These edges are supplied by weights, taken i.i.d.~with probability density $f$,
whereas at the same time, the weights to all other edges emanating from $i_0,i_1,\ldots,i_{k-1}$ are updated by subtracting $T_k$.
We assign marks analogously as before by drawing them uniformly at random from the set $[n]\setminus \{i_k\}$. Obviously, this
mark set is denoted by $I_k$.

\smallskip

We thin the MBPT defined above as follows. At each time point $A_k=T_1+T_2+\ldots+T_k$ we delete the newly found
individual with mark $i_k$ and the entire tree emanating from $i_k$, when
    \begin{equation}
    \label{thin-cond}
    i_k\in \{i_0,i_1,\ldots,i_{k-1}\},
    \end{equation}
i.e., we delete $i_k$ and all its offspring when the mark $i_k$ appeared previously. Observe that the probability
of the event \eqref{thin-cond} is at most $k/n$.
The sets $\{{\hat I}_k\}_{k\ge 0}$ are obtained from $\{I_k\}_{k\ge 0}$, by deleting the thinned vertices.

Obviously, for each fixed $n$, the thinned MBPT will become empty after a finite time, even in the super critical case ($(n-1)p_n>1$), when the 
{\it unthinned} MBPT will survive with positive probability. 
This is because we sample the marks from the finite set $[n]$, 
and hence with probability one at a certain random time all marks will have appeared. We consider the MBPT and its thinned version up to
the first generation that the thinned process becomes empty.
Furthermore, consider the ERRG with $n$ vertices and attach weights to the edges independently (also
independent from the branching process) with the same marginal density $f$ as used in the branching process.
Then, by the independence of binomial random variables $X_1,X_2,\ldots$ and the thinning described above,
for each $t\ge 0$, the set of marks that are reached by time $t$ in the thinned
MBPT and the set of vertices that are reached by time $t$ in the ERRG, are equal in distribution. Since
this holds for each $t\ge 0$ the proof of the following lemma is obvious:

\begin{Lemma}[FPP on ERRG is thinned CTMBP]
\label{eq-distr}
Fix $n\ge 1$ and $p\in (0,1)$.
Consider the thinned MBPT until extinction with i.i.d.~weights having density $f$ with support
$(0,\infty)$. If we apply independently an identical weight construction to the ERRG, then
for any $i_0\in [n]$, the weight $W_n(i_0,j)$ of the shortest path between initial vertex
$i_0$ and any other vertex $j\in [n]$ in the ERRG is
equal in distribution to the weight of the shortest path between
$i_0$ and $j$ in the thinned MBPT.
\end{Lemma}

\begin{Remark}
\label{coupl-hopcount}
(a) In both random environments (thinned MBPT and ERRG) it can happen that $i_0$ and $j$ are not connected.
We then put the weight equal to $+\infty$.\\
(b) Observe that since the weight distribution admits a density and is defined on finite objects,
the minimal weight, if finite, uniquely identifies the minimal path and hence the number of edges
on this minimal path, the hopcount. Hence, Lemma \ref{eq-distr} also proves that
the hopcount between $i_0$ and $j$ in the ERRG is
equal in distribution to the hopcount between
$i_0$ and $j$ in the thinned MBPT.
\end{Remark}

\subsection{The continuous-time branching process}
\label{CTBP} We study distance between uniformly chosen, {\it
connected}, vertices in $G_n(p_n)$ with exponentially distributed
edge-weights. Since the vertices of the ERRG are exchangeable we
can as well study the distance between vertex $1$ and vertex $2$
conditioned on the event that $1$ and $2$ are in the giant component. For the
remainder of the paper we put: \eqn{ \label{exp-density}
f(x)=\e^{-x},\qquad x>0, } the density of the exponential
distribution with rate $1$. For the initial vertex we take either
$i_0=1$ or $i_0=2$, and, if necessary, we will use a superscript
to indicate from which vertex we start. Since we have a ({\it
distributional}) embedding of the vertices of $G_n(p)$, reached by
time $t$,  in the marked branching process tree (see  Lemma
\ref{eq-distr} and Remark \ref{coupl-hopcount}) we can restrict
our attention exclusively to the MBPT and its thinned version.

We start from the individual with mark $1$ and let fluid percolate at unit speed from this
individual until it reaches one of the children. Because the weights are exponentially distributed (see (\ref{exp-density})),
the number of individuals that has been reached by time $t\ge 0$ form
a {\it continuous time} branching process (CTMBP) \cite[Chapter 2]{athreya}.
The {\it alive} individuals are defined as those individuals which are directly connected to the reached individuals
appropriately called {\it wetted} individuals. At time $t=0$ this number is
$S_1=X_1$. Then at the first splitting time $T_1$
at which time the second individual is wetted the number of alive vertices equals
$S_2=X_1+X_2-1$, since $X_2$ new individuals are born and one individual gave birth and died (became wet).
So, starting from the i.i.d.~sequence $X_1,X_2,\ldots,$ we define the random walk
\eqn{
\label{random-walk}
S_k=X_1+\ldots+X_k-(k-1),\quad k\ge 1
}
and the inter-splitting times $T_k,\,k\ge 1$, where conditioned on $X_1,X_2,\ldots,X_k$, the random variables $T_1, T_2,\ldots,T_k$
are independent and where $T_j$ has an exponential distribution with parameter $S_j$.
Indeed, by the memoryless property of the exponential distribution, when $t=A_j=T_1+\ldots+T_j$, we have $S_{j+1}$ {\it alive} individuals and all their $S_{j+1}$ edges, with which they are attached to the dead or wet individuals have an independent exponentially
distributed random variable with rate $1$ as weight. Therefore the inter-splitting time $T_{j+1}$
until the next individual is wetted, conditioned on $S_{j+1}$, is equal in distribution to
$\min_{1\leq s\leq S_{j+1}} E_s,
$
where $E_1,E_2,\ldots$ are independent rate $1$ exponentially
distributed random variables. This shows that the conditional distribution of
$T_{j+1}$ is the correct one, since the minimum of $q$ independent $\EXP(1)$ has an exponential
distribution with rate $q$.

At each splitting time $A_k, \, k\ge 1$, the number of wetted individuals equals $k+1$, whereas the
number of edges and endpoints competing to become wet equals $S_{k+1}$.
We now introduce {\it shortest-weight trees}. Since the ERRG is imbedded in the CTMBP, we introduce the
$\BP$ in terms of marks and splitting times. The shortest weight tree $\BP_k$ is equal to
the collection of marks that are wetted at the $k^{\rm th}$ splitting time $A_k$ and we include 
in the definition the splitting-times $A_1,A_2,\ldots,A_k$. So, $\BP_0=(\{i_0\},A_0=0)$, the mark of the root,
and
\eqn{
\label{def-BP}
\BP_k=\Big(\{i_0,i_1,\ldots,i_k\},\{A_0,A_1,\ldots,A_k\}\Big),\quad k\ge 1,
}
where $i_0,i_1,\ldots,i_k$, denote the marks reached at the splitting times $A_0,A_1,\ldots,A_k$, respectively.
To connect $i_j$ with the splitting time $A_j$, we introduce the mapping $t(i_j)=j$, so that mark $i\in \BP_k$ was reached at time
$A_{t(i)}$.
For a fixed mark $i\in\BP_k$ and all marks $j\neq i$,
the random variables $\ind_{ij}|\BP_k$ that indicate whether edge $ij$ is in the branching process tree
(or in other words: $j$ is one of the alive individuals born out of $i$) are \emph{independent} and satisfy
    \eqn{
    \label{succesprob}
    \prob\Big(\ind_{ij}=1|\BP_k\Big)=p_n.
    }

Furthermore, let us fix $k\ge 0$ and consider the in total $S_{k+1}$ alive individuals at the splitting time $A_k$.
The edges between the wet mark $i$ and the mark $j$, corresponding to an alive individual with mark $j$, has weight
$E_{ij}$, $i\in \BP_k$. Now observe from the memoryless property of the exponential distribution, that,
conditionally on $\BP_k$,
    \eqn{
    \label{memoryless-prop}
    E_{ij}| \BP_k\stackrel{d}{=}A_{k+1}-A_{t(i)}+E'_{ij},
    }
where $E'_{ij}$ are independent $\EXP(1)$ random variables.

\subsection{Hopcount and weight from root to newly added individuals}
\label{HWRNV}
Concerning the height and the weight of the path with minimal weight from root $1$
in  $\{\BP_k\}_{k\ge 0}$ to the $m^{\rm th}$ individual that has been wetted we refer to
\eqref{Gm-def} and \eqref{Tm-def}, respectively, where the random variables $I_1,I_2,\ldots,I_m$
are conditionally independent given $X_1,X_2
, \ldots, X_m$, the binomially distributed offspring in the tree with root $1$,
 and where $d_j=X_j$, $s_j=S_j$. Indeed, the height of the $m^{\rm th}$ individual that has
 been wetted in the CTMBP is equal to the generation $G_m$ of that individual
and its weight is equal to the splitting-time $A_m=T_1+\ldots+T_m,$
where conditionally on $X_1,X_2, \ldots, X_m$, the
distribution of $T_j$ is given by
    \eqn{
        T_j \stackrel{d}{=} E_j/S_j,
    }
and where $\{E_i\}_{i=1}^{\infty}$ are i.i.d.\ exponential random
variables with mean 1, independent of all random variables introduced earlier.
The CLT for the hopcount $G_m$ and the limit law for the weight $A_m$ in the CTMBP can be deduced from
Theorem \ref{theo:poisson} in the following way:

\begin{Proposition}[Asymptotics for shortest weight paths in the CTMBP]
\label{propnorm}
Consider\\ an i.i.d.~sequence $X_1,X_2,\ldots,$ with
a binomial distribution where the parameters $n-1$ and $p_n$ satisfy $\lim_{n\to \infty} np_n=\lambda>1$, and put
$\beta=\lambda/(\lambda -1)$. Further more let $G_m$ and $A_m$ be defined as in
$\eqref{Gm-def}$ and $\eqref{Tm-def}$, respectively, with $d_i=X_i, i\ge 1$. Then,
conditionally on $S_i=X_1+\ldots+X_i-(i-1) > 0$ for all $1\leq i\leq m$, the asymptotics in Theorem
\ref{theo:poisson} (a-b) remain to hold, where the limits are independent.
\end{Proposition}

{\bf Proof.}
According to  \cite[Theorem 2.9]{Hofs08}, we can couple $X\sim \Bin(n-1,p_n)$ and $D\sim \Poi((n-1)p_n)$,
by means of a pair of variables $({\tilde X},{\tilde D})$ with ${\tilde X}\sim \Bin(n-1,p_n)$ and ${\tilde D}\sim \Poi((n-1)p_n)$, so that
for any $\varepsilon>0$, and $n$ sufficiently large,
\eqn{
\prob({\tilde X}\neq {\tilde D})\leq (n-1)p_n^2\leq \frac{(\lambda+\varepsilon)^2}{n}.
}
Hence for $m=o(n)$, we can couple the i.i.d.~sequence $X_1,X_2,\ldots,X_m$ to an i.i.d.~sequence of Poisson variables
$D_1,D_2,\ldots,D_m$, such that
\eqn{
\lbeq{couplingofbintopoiss}
\prob({\tilde X_i}\neq {\tilde D_i},\,\mbox{for some $i\le m$})\leq \sum_{i=1}^m\frac{(\lambda+\varepsilon)^2}{n}=m\frac{(\lambda+\varepsilon)^2}{n}
\to 0.
}
Next, we define $(J_i,K_i),1\le i\le m$, conditioned on $({\tilde X_1},{\tilde D_1}),\ldots,({\tilde X_m},{\tilde D_m})$ as follows.
Let $U_1,U_2,\ldots,U_m$ be an i.i.d.~sequence of uniform $(0,1)$ random variables.
For each $i$, the sample space of $(J_i,K_i)$ is $\{0,1\}^2$, and the conditional probabilities are defined by
\eqn{
\label{joint2by2}
\prob(J_i=1\mid U_i)={\bf 1}_{\{U_i\le \frac{{\hat X}_i}{{\hat S}_i^{\sss X}}\}}, \quad
\prob(K_i=1\mid U_i)={\bf 1}_{\{U_i\le \frac{{\hat D}_i}{{\hat S}_i^{\sss D}}\}}
}
where ${\hat S}_i^{\sss X}=({\hat X}_1+\ldots+{\hat X}_i)-(i-1)$, ${\hat S}_i^{\sss D}=({\hat D}_1+\ldots+{\hat D}_i)-(i-1)$,
and where ${\bf 1}_{\cal A}$ denotes the indicator of the set ${\cal A}$.
Note that the joint distribution of $(J_i,K_i)$, conditioned on $({\tilde X_1},{\tilde D_1}),\ldots,({\tilde X_m},{\tilde D}_m)$,
is completely specified by the probabilities in \eqref{joint2by2}. Finally, set
    \eqn{
    \label{annn}
    {\cal A}_m=\{({\tilde X_1}={\tilde D_1}),\ldots,({\tilde X_m}={\tilde D_m})\}
    }
Then by \eqref{couplingofbintopoiss}, we have that $\prob({\cal A}_{m_n})\to 1$, and
    \begin{eqnarray}
    \label{distr-eq}
    &&\prob\Big(\sum_{i=1}^{m_n} J_i=\sum_{i=1}^{m_n} K_i\Big)\ge\prob({\cal A}_{m_n})\to 1,
    \end{eqnarray}
because on ${\cal A}_{m_n}$ we have $J_i=K_i, \, 1\le i \le m_n$.
In order to prove Part (a), we have to show that $\sum_{i=1}^{m_n} J_i$ conditioned on ${\hat S}_i^{\sss X}>0, \, 1\le i\le m_n$ has, {\bf whp}, the
same distribution as $\sum_{i=1}^{m_n} K_i$ conditioned on ${\hat S}_i^{\sss D}>0, \, 1\le i\le m_n$. Observe that the statement
follows from \eqref{annn} and \eqref{distr-eq} since by the coupling introduced above,
    \begin{eqnarray}
    &\prob({\hat S_i}^{\sss D}>0, \, 1\le i\le m_n)\to 1-p_{\lambda}>0.
    \label{cond2}
    \end{eqnarray}
Hence, Part (a) follows  from Part (a) of Theorem  \ref{theo:poisson}.
Referring to Part (b) of Theorem  \ref{theo:poisson} and the definition of $A_m$ in \eqref{Tm-def}, we conclude that the
proof of Part (b) also follows in a straightforward manner. \hfill\qed

\subsection{The connection time}
\label{TCT}
In Definition \eqref{def-BP} of Section \ref{CTBP} we have introduced
$\{\BP_k\}_{k\ge 0}$, which includes the set of marks of the individuals that
are wet after the $k^{\rm th}$ split.
More precisely, we grow a  CTMBP  from a root with mark $1$ and include
the marks that are successively reached and their splitting-times in $\{\BP^{\sss (1)}_k\}_{k\ge 0}$.
Our plan is to grow $\BP^{\sss (1)}_k$ until the set of marks has reached size  $a_n=\lceil \sqrt{n}\rceil$;
the branching process then contains $a_n$ wet individuals and the same marks can appear more than once.
Then, we grow independently, starting from a root with mark $2$, a second marked branching process.
We denote the shortest weight tree, i.e., the wetted marks and their splitting times of this second process by $\{\BP^{\sss (2)}_k\}_{k\ge 0}$.
We stop the first time that this second process contains a mark from $\BP^{\sss (1)}_{a_n}$, i.e.,
we stop growing $\{\BP^{\sss (2)}_k\}_{k\ge 0}$ at the {\it random} time $\CE_n$ defined by
    \eqn{
    \label{conn-time}
    \CE_n=\min \{m\ge 0: \BP^{\sss (2)}_m \cap \BP^{\sss (1)}_{a_n}\neq \varnothing\},
    }
where $\BP^{\sss (2)}_m \cap \BP^{\sss (1)}_{a_n}$ denotes the common marks in $\BP^{\sss (2)}_m$ and $\BP^{\sss (1)}_{a_n}$.
We now reach the {\it main} theorem of this section, where we establish the connection
between the weight $W_n$ and the hopcount $H_n$ in the ERRG  and the
weight and height in the shortest-weight trees $\{\BP^{\sss (i)}_k\}_{k\ge 0},\, i=1,2,$
introduced in this section.

\begin{Theorem}[Connecting hopcount to heights in CTMBP]
\label{main-theor}
Let $H_n$, $W_n$ denote the hopcount and the weight, respectively, of the minimal path
between the vertices $1$ and $2$ in the Erd\H{o}s-R\'enyi graph $G_n(p_n)$, where we condition on the vertices $1$ and $2$ to be in the giant component.\\
(a) For $n\to \infty$, {\bf whp},
\eqn{
\label{hopid}
H_n\stackrel{d}{=} G^{\sss (1)}_{U_n,a_n}+G^{\sss (2)}_{\CE_n},
}
where $G_{U_n,a_n}$ is the height (or generation) of a uniformly chosen mark in $\BP^{\sss (1)}_{a_n}$,
and $G^{\sss (2)}_{\CE_n}$ is the height of the mark attached at the random time
$\CE_n$  defined in \eqref{conn-time}, and where conditioned on $\CE_n$ the variables $G^{\sss (1)}_{U_n,a_n}$ and $G^{\sss (2)}_{\CE_n}$
are independent. \\
(b) Similarly, for $n\to \infty$, {\bf whp},
\eqn{
\label{weightid}
W_n\stackrel{d}{=} A^{\sss (1)}_{a_n}+A^{\sss (2)}_{\CE_n},
}
where $A^{\sss (i)}_{m},\,i=1,2,$  is the splitting-time of the $m^{\rm th}$ mark in  $\{\BP^{\sss (i)}_{k}\}_{k\ge 0}$.
\end{Theorem}

For the moment, we postpone the proof of Theorem \ref{main-theor}. We first need to establish two intermediate results.
Observe that if in $\BP^{\sss (i)}_k,\, i=1,2,$ all marks on the minimal weight path between the root $i$ and
some other mark $U\in \BP^{\sss (i)}_k$ appeared for the first time, then the entire path is contained in the thinned branching process, and
by Theorem  \ref{eq-distr} and the remark following that same theorem both the weight and the
hopcount in the ERRG between the vertices $i$ and $U$ are equal in distribution
to $G^{\sss (i)}_{U}$ and $A^{\sss (i)}_{U}$, where $G^{\sss (i)}_{U}$ denotes the generation of
mark $U$ and $A^{\sss (i)}_{U}$ denotes the splitting-time of $U$. Hence, in order to prove
Theorem \ref{main-theor}, we need an upper bound on the expected number
of marks in $\BP^{\sss (i)}_k$ whose minimal weight path contains thinned marks. This will be the content
of the first lemma below. The second lemma below shows that $\CE_n/a_n$ converges to an $\EXP(1)$ random variable, so that we have a handle
on the size of $\CE_n$.

\begin{Lemma}[A coupling bound]
Fix $k\ge 1$ and denote by $M_k^{\sss (i)},\,i=1,2,$ the number of marks in $\BP^{\sss (i)}_k$ for which the minimal path
from root $i$ to that mark contains a thinned mark (or individual). Then
    \eqn{
    \label{mult-marks}
    \expec[M_k^{\sss (i)}]\le \frac{k^2}{n-a_n},\qquad i=1,2.
    }
The above inequality holds in $\BP^{\sss (1)}_k$ for $1\le k\le a_n$ and  in
$\BP^{\sss (2)}_k$ for $1\le k \le \CE_n$.
\end{Lemma}

\noindent{\bf Proof.}
We start with $\BP^{\sss (1)}_k$. Let $\BP^{\sss (1)}_k$ contain $M_k^{\sss (1)}$ marks
whose minimal path contains at least one thinned mark. The new mark which is drawn to be the mark of the newly wetted
vertex in $\BP^{\sss (1)}_{k+1}$ attaches to one of these $M_k^{\sss (1)}$ marks with probability $M_k^{\sss (1)}/(k+1)$, since $|\BP^{\sss (1)}_k|=k+1$.
If the new mark is attached to a vertex with one of the $k+1-M_k^{\sss (1)}$ other marks, then with probability at most
$(k+1)/n$ this mark appeared previously and has to be thinned. Hence
    \eqn{
    \label{recursion-Nk}
    \expec[M^{\sss (1)}_{k+1}-M^{\sss (1)}_k|M^{\sss (1)}_k]\le \frac{k+1}{n}+\frac{M^{\sss (1)}_k}{k+1}.
    }
For $\BP^{\sss (2)}_k$ the fraction $\frac{k+1}{n}$ in the recursion should be replaced by $\frac{k+1}{n-a_n}$, because
for each $1\le k\le\CE_n$ we know that $\BP^{\sss (1)}_{a_n}$ contains at most $a_n$ different marks.
Write $b_k=\expec[M_k^{\sss (1)}]$. Taking double expectations in \eqref{recursion-Nk} and solving the recursive inequality yields:
    $$
    \frac{b_{k+1}}{k+1}\leq \frac{1}{n}+\frac{b_k}{k}, \qquad \mbox{so that} \qquad \frac{b_{k+1}}{k+1}\le \frac{k}{n}.
    $$
This yields $\expec[M_k^{\sss (1)}]\le \frac{k^2}{n}$. The proof for $\BP^{\sss (2)}_k$ is similar, however since at most
$a_n$ marks are occupied by $\BP^{\sss (1)}_{a_n}$, we replace $n$ by $n-a_n$ and this yields \eqref{mult-marks}.
\hfill\qed

\begin{Lemma}[Weak convergence of connection time]
Conditioned on the event that both CTMBPs survive, the connection time
$\CE_n$, defined in \eqref{conn-time}, satisfies the asymptotics
\eqn{
\label{expconn}
\CE_n/a_n \stackrel{d}{\to} E,
}
where $E$ has an exponential distribution with rate 1.
\end{Lemma}

\noindent {\bf Proof.}
Define by $\Qprob^{\sss(j)}_n$ the conditional probability given
both $\BP^{\sss (1)}_{a_n}$ and $\BP^{\sss (2)}_j$. Similar to \cite[Lemma B.1]{BHHS08}, the probability
$\prob(\CE_n>m)$ satisfies the following product form:
 \eqn{
    \lbeq{prod-CEn}
    \prob(\CE_n>m)=\expec\Big[\prod_{j=1}^m \Qprob^{\sss(j)}_n(\CE_n>j|\CE_n>j-1)\Big].
    }
We omit the proof of \eqref{prod-CEn}, which follows by suitable conditioning arguments and is identical to
\cite[Proof of Lemma B.1]{BHHS08}. Since both branching processes survive, $\BP^{\sss (1)}_{a_n}$ contains
$a_n(1+o(1))$ different marks and by \eqref{mult-marks} the expected number of multiple marks is $O(1)$. From this we obtain that
as long as the second branching process does not die out,
    $$
    \Qprob^{\sss(j)}_n(\CE_n>j|\CE_n>j-1)=1-\frac{a_n(1+o(1))}{n}.
    $$
Substitution of this into \refeq{prod-CEn} with $m=a_n x$, and $x>0$, yields
    \eqn{
    \prob(\CE_n> a_n x)=\Big[1-\frac{a_n(1-o(1))}{n}\Big]^{a_n x}\to \e^{-x},
        }
and we arrive at \eqref{expconn}.\hfill\qed

\paragraph{Proof of Theorem \ref{main-theor}.}
We start with the proof of the hopcount.
Since $1$ and $2$ are connected, with high probability,  $1$ and $2$ are contained in the
giant component of the ERRG. This implies that both  CTMBP's
do not die out. Consider the mark set $\BP^{\sss (1)}_{a_n}$ and
choose one of the wet individuals at random. The mark number of this individual will be denoted by
$U_n$. Then the probability that the shortest weight path from root $1$ to $U_n$ contains a thinned mark is at most
    \eqn{
    \frac{\expec[M^{\sss  (1)}_{a_n}]}{a_n}\le \frac{a_n}{n-a_n}\to 0.
    }
This shows that {\bf whp} the hopcount between vertex
$1$ and vertex $U_n$ in $G_n(p_n)$ is given by
$G^{\sss (1)}_{U_n,a_n}$. By \eqref{expconn}, $\CE_n=\op(\sqrt{n}\log n),$
and a similar argument shows that, {\bf whp}, the hopcount
between $2$ and the individual corresponding to the last mark added to
$\BP^{\sss (2)}_{\CE_n}$ is equal to
$G^{\sss (2)}_{\CE_n}$.

Now recall the discussion around \eqref{succesprob}
and \eqref{memoryless-prop}. We conclude from this discussion that
indeed at time $\CE_n$ one of the marks in $\BP^{\sss (1)}_{a_n}$ is chosen
uniformly at random and this establishes \eqref{hopid}. 
Finally, we  introduce the notation $p(i)$ for $i\in \BP_k$, to denote the parent of mark $i$,
i.e., the mark to which $i$ was attached in $\BP_k$. Because of the way the exponential weight in  \eqref{memoryless-prop}
can be re-allocated to edges in $\BP^{\sss (1)}_{a_n}$, the edge $\{p(i),i\}$
gets weight $A_{a_n+1}-A_{t(i)}$ and
the weight $E'_{ij}$ is allocated to $ij$, with $i\in \BP^{\sss (1)}_{a_n}$ and
$j\notin \BP^{\sss (1)}_{a_n}$, so that we get
\eqref{weightid}. \hfill \qed

\section{Proof of the main results}
\label{sec-proof-main-res}
In this section, we complete the proof of our main results.

\noindent{\bf Proof of Theorem \ref{theo:hopcount}.}
Observe that the statistical dependence between the two CTMBPs, i.e., the one that grows  from
the root with mark 1 and the one that grows from the root with mark 2, is introduced {\it only} through their marks. Moreover, in
each tree, the hopcount of the $m^{\rm th}$ added individual is independent of the marks.
Hence, conditioned on the event $\{\CE_n=m\},$ the random variables $G^{\sss (1)}_{U_n,a_n}$ and
$G^{\sss (2)}_{\CE_n}=G^{\sss (2)}_{m}$ are statistically independent.

Denote by $Z^{\sss (1)}_n=(G^{\sss (1)}_{U_n,a_n}-\beta\log{a_n})/\sqrt{\beta \log{a_n}}$
and by $Z^{\sss (2)}_n=(G^{\sss (2)}_{\CE_n}-\beta\log{\CE_n})/\sqrt{\beta \log{\CE_n}}$,
and by $\xi_1$, $\xi_2$ two independent standard normal
random variables. Then, by Proposition \ref{propnorm}, for any sequence $m_n=o(n)\to \infty$,
and for any bounded continuous function $g$ of two real variables
\eqn{
\lim_{n\to \infty}\expec[g(Z^{\sss (1)}_n,Z^{\sss (2)}_n)|\CE_n=m_n]=\expec[g(\xi_1,\xi_2)].
}
Hence by bounded convergence (the function $g$ is bounded), and using \eqref{expconn},
\begin{eqnarray*}
\lim_{n\to \infty}\expec[g(Z^{\sss (1)}_n,Z^{\sss (2)}_n)]
&=& \int_0^\infty
\lim_{n\to \infty}\expec\Big[g(Z^{\sss (1)}_n,Z^{\sss (2)}_n)|\CE_n=\lceil a_n y \rceil\Big] \,d\prob(\CE_n\leq \lceil a_n y \rceil)\\
&=&\int_0^\infty \expec[g(\xi_1,\xi_2)]\e^{-y}\,dy=\expec[g(\xi_1,\xi_2)].
\end{eqnarray*}
This shows the joint weak convergence of the pair $(Z^{\sss (1)}_n,Z^{\sss (2)}_n)$ to $(\xi_1,\xi_2)$.
By the continuous mapping theorem \cite[Theorem 5.1]{Bill99} applied to $(x,y)\mapsto x+y$,
we then find that
\begin{eqnarray}
\frac{H_n-\beta \log n}{\sqrt{\log n}}
&=& \frac{G^{\sss (1)}_{U_n,a_n}-(\beta/2) \log n}{\sqrt{\log n}}+
\frac{G^{\sss (2)}_{\CE_n}-(\beta/2) \log n}{\sqrt{\log n}}\nonumber\\
&=& \frac{G^{\sss (1)}_{U_n,a_n}-\beta \log a_n}{\sqrt{\log a_n}}\cdot \frac{\sqrt{\log a_n}}{\sqrt{\log n}}
+\frac{G^{\sss (2)}_{\CE_n}-\beta \log \CE_n}{\sqrt{\log \CE_n}}
\cdot \frac{\sqrt{\log \CE_n}}{\sqrt{\log n}}+\op(1)\nonumber\\
&\stackrel{d}{\to}& \xi_1/\sqrt{2}+\xi_2/\sqrt{2},
\end{eqnarray}
and we note that $(\xi_1+\xi_2)/\sqrt{2}$ is again standard normal. \hfill \qed

\noindent{\bf Proof of Theorem \ref{theo:weight}.}
%Conditioned on $S_1>0,\ldots,S_m>0$, we obtain, as in the proof of
%Proposition \ref{propnorm}, that
%    \eqn{
%    A_m-\gamma \log m \convas X=-\gamma
%    \log (\gamma W_{\lambda}),
%    }
%where $W_{\lambda}$ was introduced in (2.3) and where $\gamma=\frac1{\lambda-1}$.
Similarly to the proof of Theorem \ref{theo:hopcount}, for any sequence $m_n\to \infty$,
\eqn{
\lbeq{joint-limit}
\Big(A^{\sss (1)}_{a_n}-\gamma \log a_n,
A^{\sss (2)}_{\CE_n}-\gamma \log \CE_n|\CE_n=m_n\Big)
\convd
(X_1,X_2),
}
with $X_1,X_2$ two independent copies of $X=-\gamma\log (\gamma W_{\lambda})$.
We rewrite the second coordinate as:
$$
A^{\sss (2)}_{\CE_n}-\gamma \log \CE_n
=
A^{\sss (2)}_{\CE_n}-\gamma \log a_n-\gamma \log(\CE_n/a_n).
$$
Hence, from \eqref{joint-limit} and the limit law for $\CE_n/a_n$ given in \eqref{expconn},
we obtain:
\begin{eqnarray*}
W_n-\gamma \log n &=&A^{\sss (1)}_{a_n}+A^{\sss (2)}_{\CE_n}-2\gamma\log a_n\\
&=&A^{\sss (1)}_{a_n}-\gamma\log a_n+A^{\sss (2)}_{\CE_n}-\gamma\log \CE_n+\gamma \log(\CE_n/a_n)\\
&\stackrel{d}{\to}& X_1+X_2+\gamma \log(E),
\end{eqnarray*}
where $E\stackrel{d}{=} \EXP(1)$ is independent of $(X_1, X_2)$. We can reformulate this as
\eqn{
(\lambda-1)W_n- \log n
\convd -\log (\gamma W^{\sss (1)}_{\lambda})-\log (\gamma W^{\sss (2)}_{\lambda})-M,
\lbeq{ult-form}
}
where $M=-\log(E)$ has a Gumbel distribution, i.e., $\prob(M\leq x)=\Lambda(x)=\exp(-\e^{-x})$. \hfill \qed

\noindent{\bf Proof of Corollary \ref{cor-lim}.}
We start with a proof of
\eqref{weakconvG}, with $\beta$ replaced by $\beta_n=\lambda_n/(\lambda_n-1)$. We give the proof for the 
CTMBP with $\Bin(n,\lambda_n/n)$ offspring distribution
(thus avoiding the coupling with Poisson random variables in Proposition \ref{propnorm}).
Since, for $\lambda_n\to \infty$, any two vertices are, {\bf whp}, connected
it is not necessary to condition on the event $\{S_i>0,\,1\le i \le m\}$, and therefore it is straightforward from \eqref{Gm-def} and
\eqref{bi-def} that
    \eqn{
    \lbeq{def-Gvann}
    G_m \stackrel{d}{=} \sum_{i=1}^m I_i,\qquad
    \mbox{where}\qquad
    \prob(I_i = 1|\{D_i^{\sss (n)}\}_{1\le i \le m}) = D_i^{\sss (n)}/S_i^{\sss (n)},
    }
with $D_1^{\sss (n)},D_2^{\sss (n)},\ldots,D_m^{\sss (n)}$, independent $\Bin(n,\lambda_n/n)$ and
$S_i^{\sss (n)}=\sum_{j=1}^i D_j^{\sss (n)} -(i-1)$. By \cite[Chapter 3]{JanLucRuc00},
the random variables $D_j^{\sss (n)}$ concentrate around $\lambda_n$, whereas
$S_i^{\sss (n)}$ concentrates around $i(\lambda_n-1)+1$.
Therefore $G_m$ in \eqref{def-Gvann} satisfies the same asympotics as
    \[
    \sum_{i=1}^m {\hat I}_i^{\sss (n)},\qquad
    \mbox{where}\qquad
    \prob({\hat I}_i^{\sss (n)} = 1) = \frac{\lambda_n}{i(\lambda_n-1)}=\frac{\beta_n}{i},
    \]
with ${\hat I}_1^{\sss (n)},{\hat I}_2^{\sss (n)},\ldots,{\hat I}_m^{\sss (n)}$ independent. Applying Lemma \ref{lemma:normal},
we conclude that for $m=m_n \to \infty$,
    \eqn{
    \lbeq{interres}
    \frac{G_m - \beta_n\sum_{i=1}^m 1/i}
    {\sqrt{\sum_{i=1}^m (\beta_n/i)(1-\beta_n/i)}}
    \convd Z,
    }
where $Z$ is standard normal. Note that the denominator of \eqref{interres} can be replaced by
$\sqrt{\log m}$, since $\beta_n\to 1$ and hence $\sum_{i=1}^m (\beta_n/i)(1-\beta_n/i)/\log m\to 1$.
This proves \eqref{weakconvG}, with $\beta$ replaced by $\beta_n=\lambda_n/(\lambda_n-1)$.
The centering constant $\beta_n\sum_{i=1}^m 1/i$ can be replaced by $\beta_n \log m$, but
in general not by $\log m$, except when $(\beta_n-1)\sqrt{\log m_n}\to 0$, or equivalently
$\lambda_n=o(\sqrt{\log m_n})$.

We now turn to the proof of Part (b) of Corollary \ref{cor-lim}. Again, we will use that
$S_i^{\sss (n)}$ concentrates around $i(\lambda_n-1)+1$. Hence, from \refeq{Tm-def},
$$
A_m\stackrel{d}{=}\sum_{i=1}^m E_i/S_i^{\sss (n)}\approx \frac{1}{\lambda_n-1}\sum_{i=1}^m E_i/i.
$$
Now define $B_n$ as a random variable with distribution equal to the maximum of $n$ independent exponentially
distributed random variables with rate $1$, i.e.,
$$
\prob(B_n \le x)=(1-\e^{-x})^n.
$$
Then $E_i/i,\, 1\le i \le n,$ are equal in distribution to the spacings of these
exponentially distributed variables i.e., $B_n\stackrel{d}{=}\sum_{i=1}^n E_i/i$.
Consequently, with $\gamma_n=1/(\lambda_n-1)$,
\begin{eqnarray}
&&\lim_{n\to \infty}
\prob (\gamma_n^{-1}A^{\sss (1)}_{a_n}-\log a_n\le x)=
\lim_{n\to \infty}\prob(B_{a_n}-\log a_n \le x)\nonumber\\
&&\qquad= \lim_{n\to \infty} (1-\e^{-x+\log a_n})^{a_n}=\Lambda(x),
\end{eqnarray}
where $\Lambda$ denotes the distribution function of a Gumbel random variable, i.e.,
$\Lambda(x)=\exp(-\e^{-x})$. Parallel to the proof of Theorem \ref{theo:weight},
we conclude that
$$
\big(\gamma_n^{-1}A^{\sss (1)}_{a_n}-\log a_n,
\gamma_n^{-1}A^{\sss (2)}_{\CE_n}- \log \CE_n
\big)
\stackrel{d}{\to}
(M_1,M_2),
$$
where $M_1$ and $M_2$ are two independent copies of a random variable with distribution
function $\Lambda$. Finally, since $\CE_n/a_n$ converges in distribution to
an $\EXP(1)$ random variable, we conclude from the
continuous mapping theorem that \eqref{weightconlamn} holds. \hfill \qed

\noindent{\bf Proof of Theorem \ref{theo:extrema}.}
This proof is a consequence of Theorems \ref{theo:hopcount}-\ref{theo:weight}
and the results on the diameter of $G_n(\lambda/n)$ proved in
\cite{FerRam04}. We shall sketch the proof.
We start by recalling some notation and properties of the
ERRG from \cite{FerRam04}. Let the \emph{2-core} of $G_n(\lambda/n)$
be the maximal subgraph of which each vertex has degree at least 2.
The study of the diameter in \cite{FerRam04} is based on the crucial
observation that any longest shortest path in a random graph
(with minimum degree 1) will occur between a pair vertices $u, v$ of degree
1. Moreover, this path will consist of three segments: a path from $u$ to the
2-core, a path through the 2-core, a path from the 2-core to $v$.
While this is used in \cite{FerRam04} only for graph distances, the
same applies to FPP on a graph. Now, when we pick two uniform vertices $i$ and $j$
(as in Theorems \ref{theo:hopcount}-\ref{theo:weight}), then the paths from
$i$ and $j$ to the 2-core are (a) unique and (b) {\bf whp} of length
$o(\log{n})$ (since the giant component with the 2-core removed is a collection of trees
of which most trees have size $o(\log{n})$). As a result,
Theorems \ref{theo:hopcount}-\ref{theo:weight} also hold when picking two \emph{uniform
vertices in the 2-core}.

Next we investigate the maximal weight and length of shortest-weight
paths in the ERRG. In \cite[Section 6.1]{FerRam04}, it is shown that
{\bf whp} there are two paths of length at least $(1/\log{(-\mu_\lambda)}-\vep)\log{n}$
in the giant component with the 2-core removed connecting two vertices $i^*$ and $j^*$
of degree 1 to two vertices $U,V$ in the 2-core. Since these paths are
\emph{unique}, they must also be the shortest-weight paths for FPP on the ERRG
between the vertices $i^*$ and $U$, and $j^*$ and $V$, respectively.
The weights along these paths is thus a sum of $(1/\log{(-\mu_\lambda)}-\vep)\log{n}$
i.i.d.~$\EXP(1)$ random variables. By exchangeability of the vertices, $U,V$ are
uniform vertices in the core. Thus, the shortest-weight path from $i^*$ to
$j^*$ has at least $\lambda/(\lambda-1)\log{n}+\op(\log(n))+(2/\log{(-\mu_\lambda)}-2\vep)\log{n}$
hops, as required, while its weight is at least
$1/(\lambda-1)\log{n}+\op(\log(n))+(2/\log{(-\mu_\lambda)}-2\vep)\log{n}$.
\hfill\qed

\paragraph{Acknowledgements.}
The work of RvdH was supported
in part by the Netherlands Organisation for Scientific Research (NWO).

\bibliographystyle{plain}

%\bibliography{bib-erdos}

\begin{thebibliography}{10}

\bibitem{athreya}
K.B. Athreya and P.E. Ney.
\newblock {\em {Branching Processes}}.
\newblock Dover Publications, 2004.

\bibitem{vcg-random-shanky}
S.~Bhamidi.
\newblock First passage percolation on locally tree like networks i: Dense
  random graphs.
\newblock {\em Journal of Mathematical Physics}, 49(12): 125218, 27, 2008.

\bibitem{BHHS08}
S.~Bhamidi, R.~van~der Hofstad, and G.~Hooghiemstra.
\newblock First passage percolation on random graphs with finite mean degrees.
\newblock {\em Preprint, to appear in Ann.\ Appl.\ Probab.}, 2010.

\bibitem{Bill99}
P.~Billingsley.
\newblock {\em Convergence of Probability Measures}.
\newblock Wiley Series in Probability and Mathematical Statistics. John Wiley
  \& Sons Inc., New York, second edition, 1999.
\newblock A Wiley-Interscience Publication.

\bibitem{buhler}
W.J.~B{\"u}hler.
\newblock {Generations and degree of relationship in supercritical Markov
  branching processes}.
\newblock {\em Probability Theory and Related Fields}, 18(2):141--152, 1971.

\bibitem{stein-chen}
L.H.Y.~Chen and Q.M.~Shao.
\newblock {Stein's method for normal approximation}.
\newblock {\em in: An introduction to SteinÕs method, ed. A.D. Barbour and
  L.H.Y. Chen, Lecture notes series IMS, Natinal University Singapore}, pages
  1--59, 2005.

\bibitem{dudley}
R.M.~Dudley.
\newblock {\em {Real Analysis and Probability}}.
\newblock Cambridge University Press, 2002.

\bibitem{FerRam04}
D.~Fernholz and V.~Ramachandran.
\newblock The diameter of sparse random graphs.
\newblock {\em Random Structures Algorithms}, 31(4):482--516, 2007.


\bibitem{ForKasGin71}
C.~M. Fortuin, P.~W. Kasteleyn, and J.~Ginibre.
\newblock Correlation inequalities on some partially ordered sets.
\newblock {\em Comm. Math. Phys.}, 22:89--103, 1971.


\bibitem{hamm-welsh}
J.M.~Hammersley and D.J.A.~Welsh.
\newblock {First-passage percolation, sub-additive process, stochastic network
  and generalized renewal theory}.
\newblock {\em Bernoulli, 1713: Bayes, 1763; Laplace, 1813. Anniversary
  Volume}, 1965.

\bibitem{Hofs08}
R.~van~der Hofstad.
\newblock Random graphs and complex networks.
\newblock In preparation, see {\tt
  http://www.win.tue.nl/$\sim$rhofstad/NotesRGCN.pdf}, 2009.


\bibitem{hofs-erdos-fpp}
R.~van~der Hofstad, G.~Hooghiemstra, and P.~Van~Mieghem.
\newblock First-passage percolation on the random graph.
\newblock {\em Probab. Engrg. Inform. Sci.}, 15(2):225--237, 2001.

\bibitem{howard}
C.D.~Howard.
\newblock {Models of first-passage percolation}.
\newblock {\em Probability on Discrete Structures}, pages 125--173, 2004.

\bibitem{JanLucRuc00}
S.~Janson, T.~{\L}uczak, and A.~Rucinski.
\newblock {\em Random graphs}.
\newblock Wiley-Interscience Series in Discrete Mathematics and Optimization.
  Wiley-Interscience, New York, 2000.

\bibitem{rossboek}
S.M.~Ross.
\newblock {\em {Stochastic Processes}}.
\newblock John Wiley and Sons, 1996.

\end{thebibliography}

\def\cprime{$'$}

\end{document}